 \documentclass[12pt]{article}

\usepackage{amstext    }
\usepackage{amsthm    }
\usepackage{a4}
\usepackage[mathscr]{eucal}
\usepackage{mathrsfs}

\usepackage{amsmath}
\usepackage{amssymb}
\usepackage{amscd}

\newtheorem{theorem}{Theorem}[section]

\newtheorem{proposition}[theorem]{Proposition}

\newtheorem{lemma}[theorem]{Lemma}
\newtheorem{remark}[theorem]{Remark}

\newcommand{\cali}[1]{\mathscr{#1}}

\newcommand{\supp}{{\rm supp}}

\newcommand{\ddc}{dd^c}

\newcommand{\dbar}{\overline\partial}

\newcommand{\id}{\mathop{\mathrm{id}}\nolimits}

\newcommand{\Tan}{{\rm Tan}}
 \newcommand{\area}{{\rm area}}

\newcommand{\Ec}{\cali{E}}
\newcommand{\Fc}{\cali{F}}

\newcommand{\Lc}{\cali{L}}

\newcommand{\Oc}{\cali{O}}

\newcommand{\B}{\mathbb{B}}
\newcommand{\C}{\mathbb{C}}
\newcommand{\N}{\mathbb{N}}

\newcommand{\R}{\mathbb{R}}
\renewcommand{\P}{\mathbb{P}}

\title{Equidistribution  for meromorphic maps with dominant topological degree}

\author{Tien-Cuong Dinh, Vi{\^e}t-Anh Nguy{\^e}n and  Tuyen Trung Truong}

\begin{document}

\maketitle

\begin{abstract}
 Let $f$ be a meromorphic self-map on a compact K\"ahler manifold whose topological degree is strictly larger than the other dynamical degrees. We show that repelling periodic points are equidistributed with respect to the equilibrium measure of $f$. We also describe the exceptional set of points whose backward orbits are not equidistributed. 
 \end{abstract}

\noindent
{\bf Classification AMS 2010:} 37F10, 32U40, 32H50 

\noindent
{\bf Keywords:}  meromorphic  self-map, periodic point, equidistribution, exceptional set, tangent current.

\section{Introduction} \label{introduction}

 Let $X$  be  a  compact K{\"a}hler  manifold of dimension
$k$ and $\omega$ a K\"ahler form on $X$  so normalized  that $\omega^k$ defines a probability measure on $X$.  
 Let $f:X\rightarrow X$ be a meromorphic map. We always assume that $f$ is
 {\it dominant}, i.e. its image contains a non-empty open subset of $X$.  {\it The iterate of order
$n$} of $f$ is defined by
 $f^n:=f\circ\cdots\circ f$, $n$ times, on a dense Zariski open set and extends to a  dominant meromorphic map on $X.$
 
 Let $d_p(f)$ (or $d_p$ if there is no possible confusion), $0\leq p\leq k$,  the 
dynamical degree   of  order $p$ of $f$. This is a bi-meromorphic invariant which  measures the  norm growth of the operators $(f^n)^*$
acting  on the Hodge cohomology group $H^{p,p}(X,\C)$ when $n$  tends
to infinity, see Section \ref{section_currents} for details. We always have $d_0(f)=1$. The last dynamical degree $d_k(f)$ is {\it the topological degree} of $f$: it is equal to the number of points in a generic fiber of $f$. We also denote it by $d_t(f)$ or simply by $d_t$. Throughout the paper, we assume that $f$ is with {\it dominant topological degree}\footnote{In some references, such a map is said to be with large topological degree; we think the word ``dominant'' is more appropriate.}  
in the sense that $d_t>d_p$ for every $0\leq p\leq k-1$. 

It is well-known that for such  a map $f,$ the  following weak limit  of probability measures
 $$\mu:=\lim_{n\to\infty} {1\over d_t^n}(f^n)^*\omega^k$$
 exists.  The probability measure $\mu$ is  called  
  {\it the  equilibrium measure}  of $f.$   
It has no mass on proper analytic  subsets of $X$, is {\it totally invariant}: 
  $d_t^{-1}f^*(\mu)=f_*(\mu)=\mu$  and  is exponentially mixing. 
 The measure $\mu$ is also the unique invariant measure with maximal entropy $\log d_t$. 
 We refer the reader to \cite{DinhSibony2,DinhSibony3,Guedj1} for details.
 
 \medskip
 
The first main result of this   article is  the following result.

\begin{theorem} \label{th_main_1}
Let $f:X\to X$ be a   meromorphic map with dominant topological degree $d_t$. Let $\mu$ be the equilibrium measure of $f$. 
Let $P_n$  (resp.  $RP_n$) denote  the set of  isolated periodic  (resp. repelling periodic) points of period $n$.
Let $Q_n$ denote either  $P_n$, $RP_n$ or their intersections with the support of $\mu$. 
Then $Q_n$ is asymptotically equidistributed with respect to $\mu$: we have 
$${1\over d_t^n} \sum_{a\in Q_n} \delta_a\to \mu \qquad \mbox{as} \quad n\to\infty,$$
where $\delta_a$ denotes the Dirac mass at $a$. In particular, we have $\# Q_n=d_t^n+o(d_t^n)$ as $n\to\infty$. 
\end{theorem}

The last assertion in the above theorem is an important point in our proof. Indeed, when $f$ admits positive dimensional analytic sets of periodic points, the classical Lefschetz formula does not allow to estimate the number of isolated periodic points.   
The upper bound $\# Q_n\leq d_t^n+o(d_t^n)$  is, in fact, obtained using a very recent theory of density of positive closed currents developed by Sibony and the first author in \cite{DinhSibony5}. 

For the rest of the proof, we need to construct enough repelling periodic points on the support of $\mu$. For this purpose, we will construct in Section \ref{section_inverse_branches} enough good inverse branches of balls for $f^n$ with controlled size, see Proposition \ref{prop_branches} below. The construction of inverse branches for holomorphic  discs  in projective manifolds can be obtained using a method developed by Briend-Duval in \cite{BriendDuval2}. Here we follow the approach developed by Dinh-Sibony in \cite{DinhSibony0} which also allows to carry out such a construction   for discs as  well as  balls in K\"ahler manifolds. Then an idea of Buff allows to obtain repelling periodic points \cite{Buff}.  The presence of indeterminacy sets for meromorphic maps is the source of several delicate technical points in the proof. For example, the obstruction to the existence of inverse branches for balls, at least in our approach, may be larger than the orbits of critical values and of indeterminacy loci. We will construct and use a positive closed $(1,1)$-current $R$ which allows to control this obstruction.

Note that when $X$ is a projective manifold, a weaker version of Theorem \ref{th_main_1} was stated in \cite{Guedj1}. The author's proof is, however, based on Lemma 3.3 therein 
whose proof is incomplete (the radius $r_\epsilon$ therein depends on the integer $l$) and the statement  is still an open question.  When $f$ is a holomorphic endomorphism of $\P^k$, the above theorem was obtained by Briend-Duval in \cite{BriendDuval1}. 
Their proof uses strongly the H\"older continuity of the dynamical Green function. For meromorphic maps, the dynamical Green function, even when it exists, is in general not continuous.
The same result for polynomial-like maps of dominant topological degree, in particular for a large family of rational maps on $\P^k$, was obtained by Sibony and the first author \cite{DinhSibony0}. For the case of dimension 1, see Brolin, Freire-Lopes-Ma\~n\'e, Lyubich and Tortrat \cite{Brolin, FreireLopesMane, Lyubich, Tortrat}.

Our construction of inverse branches of balls also allows to study the equidistribution of preimages of points by $f^n$. Let $I'$ denote {\it the second indeterminacy set} of $f$, i.e. the set of points $z$ such that $f^{-1}(z)$ has positive dimension. It is an analytic set of codimension at least equal to 2. The Zariski open set $X\setminus  I'$ is the set of points $a$ such that the fiber $f^{-1}(a)$ contains exactly $d_t$ points counted with multiplicity, see Section \ref{section_currents} for the definition  of the action of $f$ and $f^{-1}$  on subsets of $X$. 

Define $I'_0:=I'$, $I'_{n+1}:=I'_0\cup f(I'_n)$ for $n\geq 0$ and $I'_\infty:=\cup_{n\geq 0} I'_n$. 
Note that the set $I'_\infty$ is characterized  by the following property: the sequence of probability measures
$$\mu^a_0:=\delta_a,\quad \mu_{n+1}^a:=d_t^{-1}f^*(\mu_n^a) \quad \text{for } n\geq 0$$
is well-defined if and only if $a\not\in I_\infty'$.  We have $\mu_n^a=d_t^{-n}(f^n)^*(\delta_a)$. So $\mu_n^a$ is the probability measure equidistributed on the fiber $f^{-n}(a)$ where the points in this fiber are counted with multiplicity. 
One has to distinguish $I_\infty'$ with the set $\cup_{n\geq 0} f^n(I')$ which is a priori smaller. 

Let $I$ be {\it the (first) indeterminacy set} of $f$. Define also $I_0:=I$, $I_{n+1}:=I_0\cup f(I_n)$ for $n\geq 0$ and $I_\infty:=\cup_{n\geq 0} I_n$. The set $I_\infty\setminus I_\infty'$ consists of points $a\not\in I_\infty'$ such that the support of $\mu_n^a$ intersects $I$ for some $n\geq 0$. We will consider $a\not\in I_\infty\cup I_\infty'$.
Here is the second main result in this paper.

\begin{theorem} \label{th_main_2}
Let $f:X\to X$ and $\mu$ be as in the statement of Theorem \ref{th_main_1}.   Then  there is   a  (possibly empty) proper analytic  set $\Ec$ of $X$ such that 
for $a\not\in I_\infty\cup I'_\infty$ we have 
 $${1\over d_t^n}  (f^n)^\ast \delta_a \to \mu \qquad \mbox{as} \quad n\to\infty$$
if and only if  $a\not\in \Ec.$ 
 \end{theorem}

When $X$ is projective, it was shown by Guedj in \cite{Guedj1} that $\Ec$ is a finite or countable union of analytic sets, see also \cite{DinhSibony3} for the case of compact K\"ahler manifolds. The above theorem was obtained for holomorphic endomorphisms of $\P^k$ in \cite{BriendDuval2,DinhSibony0, FornaessSibony}. It also holds for polynomial-like maps with dominant topological degree. For the case of dimension 1, see \cite{Brolin, FreireLopesMane, Lyubich, Tortrat}.  Note that there are many meromorphic maps with $I'=\varnothing$ which are not holomorphic. For example, if $g:\widehat\P^k\to \P^k$ is a blow-up of $\P^k$ and $\pi:\widehat \P^k\to \P^k$ is a finite holomorphic map, then $\pi\circ g^{-1}$ is not holomorphic but its second indeterminacy set is empty.  For holomorphic maps  on $\P^k$, we have $I=I'=\varnothing.$

The  article  is  organized  as  follows. In Section \ref{section_currents} we prepare  the background
and fix  some terminology as  well as  recall auxiliary results concerning the actions of meromorphic maps  on currents and on cohomology and  the theory of density for positive closed currents.
Section  \ref{section_inverse_branches} is devoted  to the  construction of good inverse  branches of iterates of $f,$
which is  one of the main tools  of our work.
Using   these  inverse branches we prove  Theorems  \ref{th_main_2} in Section  \ref{section_equi_preimages}.  
 After establishing an upper bound on the number of isolated periodic points, 
 we   prove  Theorem  \ref{th_main_1} in  Section \ref{section_equi_periodic_points}. In the same section, we also explain how to obtain a lower bound for Lyapounov exponents of $\mu$ 
 from our construction of inverse branches for balls.

\medskip

\noindent
{\bf Acknowledgement. } The paper was partially prepared 
during the visit of the second author  at the University of Cologne upon  a Humboldt foundation research grant, 
and 
during the visit of the third author at the University of  Paris 6  and at the University of Paris 11 (Orsay). They would like to express their gratitude to these organizations for hospitality and  for  financial support.   
The second author also would like  to thank  Professor George Marinescu for  kind  help.


\section{Meromorphic maps and  currents}\label{section_currents}

 In this section  we define various operations  for   meromorphic maps and  positive closed currents
 on compact K\"ahler manifolds.   We also recall  some elements of the  theory of density of positive closed  currents   and  establish a preparatory result. 
We refer the  reader to  Demailly  \cite{Demailly2}, Dinh-Sibony  \cite{DinhSibony4, Sibony} and Voisin  \cite{Voisin} for basic notions  on positive closed  currents and  quasi-plurisubharmonic (quasi-p.s.h. for short) functions  and  basic  facts on K\"{a}hler geometry.
 
Let $X$ be a compact K\"ahler manifold of dimension $k$ and $\omega$ a K\"ahler form on $X$ as above.
If $T$ is a current on $X$ and $\varphi$ is a test form of right degree, the pairing $\langle T,\varphi\rangle$ denotes the value of $T$ at $\varphi$. 
 If $T$ is a positive $(p,p)$-current on $X$, its mass is given by the formula
$$\|T\|:=\langle T,\omega^{k-p}\rangle.$$
Note that when $T$ is, moreover, closed, its mass depends only on its cohomology class $\{T\}$ in $H^{p,p}(X,\R)$. Here  $H^{p,q}(X,\C)$ denotes the Hodge cohomology group of bidegree $(p,q)$ of $X$ and $H^{p,p}(X,\R):=H^{p,p}(X,\C)\cap H^{2p}(X,\R)$. 

We will write $T\leq T'$ and $T'\geq T$ for two real $(p,p)$-currents $T,T'$ if $T'-T$ is a positive current. We also write $c\leq c'$ and $c'\geq c$ for $c,c'\in H^{p,p}(X,\R)$ when $c'-c$ is the class of a positive closed $(p,p)$-current. If $V$ is an analytic subset of pure dimension $k-p$ in $X$, denote by $[V]$ the positive closed current of integration on $V$ and $\{V\}$ its cohomology class in $H^{p,p}(X,\R)$. 
 The cup-product in $H^*(X,\C)$ is  denoted by $\smallsmile.$

Consider now a dominant meromorphic map $f:X\to X$. Recall that $f$ is holomorphic on a Zariski open set and the closure $\Gamma$ of its graph in $X\times X$ is an irreducible analytic subset of dimension $k$.  Let $\pi_1$ and $\pi_2$ denote the canonical projections from
$X\times X$ onto its factors. If $A$ is a subset of $X$ define
$f(A):=\pi_2(\pi_1^{-1}(A)\cap\Gamma)$ and 
$f^{-1}(A):=\pi_1(\pi_2^{-1}(A)\cap\Gamma)$.
{\it The (first) indeterminacy set} $I$ of $f$ is the complement of the set of all points $z\in X$ such that $f(z)$ is of dimension $0,$ or equivalently, that $f(z)$ contains only one point. {\it The second indeterminacy set} $I'$ of $f$ is the complement of the set of all points $z$ such that $f^{-1}(z)$ is of dimension $0,$ or equivalently, that $f^{-1}(z)$ contains exactly $d_t$ points counted with multiplicity.
Both $I$ and $I'$ are analytic subsets of $X$ of codimension at least equal to $2.$

The map $f$ induces linear operators on forms and currents. The presence of indeterminacy locus makes these operators more delicate to handle. 
If $\varphi$ is a smooth $(p,q)$-form on $X$, then
$f^*(\varphi)$ is the $(p,q)$-current defined by
$$f^*(\varphi):=(\pi_1)_*(\pi_2^*(\varphi)\wedge [\Gamma]).$$
It is not difficult to see that $f^*(\varphi)$ is an $L^1$-form smooth
outside $I$. Its singularities along $I$ do not allow to iterate the operation in the same way. Nevertheless, the operation commutes with $\partial$ and $\dbar$. In particular, when $\varphi$ is closed or exact, so is $f^*(\varphi)$. 
 Therefore, $f^*$ induces a linear operator on $H^{p,q}(X,\C)$. We can iterate the later operator but in general we do not have $(f^*)^n= (f^n)^*$.

In the same way, the $(p,q)$-current 
$f_*(\varphi)$ is defined by 
$$f_*(\varphi):=(\pi_2)_*(\pi_1^*(\varphi)\wedge [\Gamma]).$$
This is an $L^1$-form smooth outside the critical values of $\pi_{2|\Gamma}$. The operator $f_*$ also commutes with $\partial,\dbar$ and  induces a linear operator 
$f_*$ on $H^{p,q}(X,\C)$.

Recall that {\it the dynamical degree of order $p$} of $f$ is defined by
\begin{eqnarray*}
d_p & = & \lim_{n\to\infty} \|(f^n)^*(\omega^p)\|^{1/n}=\lim_{n\to\infty} \|(f^n)_*(\omega^{k-p})\|^{1/n}\\
& = & \lim_{n\to\infty} \|(f^n)^*\|_{H^{p,p}(X,\C)}^{1/n}=\lim_{n\to\infty} \|(f^n)_*\|_{H^{k-p,k-p}(X,\C)}^{1/n}.
\end{eqnarray*}
The above limits exist and do not depend on the choice of $\omega$ nor on the norm fixed for $H^*(X,\C)$. 
They are bi-meromorphic invariants and play a central role in complex dynamics, see \cite{DinhSibony2} for details.
Recall also that a mixed version of the Hodge-Riemann theorem \cite{DinhNguyen1, Gromov1,Khovanskii,Teissier, Timorin} implies that $p\mapsto \log d_p$ is concave, i.e. $d_p^2\geq d_{p-1}d_{p+1}$. So $f$ has dominant topological degree if and only if $d_t>d_{k-1}$. 

We will now consider two particular cases of the pull-back operator $f^*$ on currents that will be used later on. If $\phi$ is a continuous function on $X$ then $f_*(\phi)$ is a bounded function on $X$ which is continuous outside $I'$. Therefore, if $\nu$ is a positive measure without mass on $I'$ we can define
$$\langle f^*(\nu),\phi\rangle := \langle \nu, f_*(\phi)\rangle.$$
It is not difficult to see that $f^*(\nu)$ is a positive measure whose mass is equal to $d_t$ times the mass of $\nu$ since $\pi_2$ restricted to $\Gamma$ defines a ramified covering of degree $d_t$   over $X\setminus I'$. 
If $\nu$ is the Dirac mass at a point $a\not\in I'$, then $f^*(\nu)$ is the sum of the Dirac masses on the fiber $f^{-1}(a)$ counted with multiplicity. If $\nu$ has no mass on $I,$ the positive measure $f_*(\nu)$ given by 
$$\langle f_*(\nu),\phi\rangle := \langle \nu, f^*(\phi)\rangle$$
  for every  continuous function $\phi$ on $X,$ 
is  well-defined and has the same mass as $\nu$. If $\nu$ is the Dirac mass at $a\not\in I$, then $f_*(\nu)$ is the Dirac mass at $f(a)$. 

The second situation concerns positive closed $(1,1)$-currents. If $T$ is such a current on $X$, we can write 
$T=\alpha+\ddc u$ where $\alpha$ is a smooth closed  $(1,1)$-form in the class $\{T\}$ and $u$ is a quasi-p.s.h. function. Then $u\circ \pi_2$ is a quasi-p.s.h. function on $\Gamma$ and we define 
$$f^*(T):=f^*(\alpha)+(\pi_1)_*(\ddc(u\circ \pi_{2|\Gamma})).$$
Using a local regularization of $T$, one can see that $f^*(T)$ is a positive closed $(1,1)$-current. The operator is linear and continuous on $T$. So using Demailly's regularization of $(1,1)$-currents on $X$ \cite{Demailly1}, we can easily check that the operator is compatible with cohomology, that is, we have $\{f^*(T)\}=f^*\{T\}$. The operator $f_*$ is defined in the same way on positive closed $(1,1)$-currents and is also compatible with the cohomology.

In the rest of this section, we recall basic facts on the theory of density of positive closed currents and give an abstract result which will allow us to bypass Lefschetz's fixed points formula in order to bound the number of periodic points. We will restrict ourselves to the simplest situation that is  needed for the present work. We will come back to this subject in a forthcoming paper. The reader is  invited   to consult \cite{DinhSibony5} for details.

Let $V$ be an irreducible  submanifold of $X$ of  dimension $l$. 
Let $\pi:\ E\to V$ denote the normal vector bundle  of $V$ in $X$.
For  a point $a\in V,$ if   $\Tan_a X$ and $\Tan_a V$ denote respectively the tangent spaces of $X$ and of $V$ at $a,$
the fiber $E_a:=\pi^{-1}(a)$  of $E$ over $a$ is canonically identified with the quotient space 
$\Tan_a X / \Tan_a V.$ The zero section of $E$  is naturally identified with $V$. 
Denote by $\overline E$ the natural compactification of $E$, i.e.
the projectivisation $\P(E\oplus \C)$ of the vector bundle $E\oplus \C$, where $\C$ is the trivial line bundle over $V$.
We still denote by $\pi$ the natural projection from $\overline E$ to $V$. 
Denote by $A_\lambda$ the multiplication by $\lambda$ on the fibers of  $E$ where $\lambda\in\C^*,$
i.e. $A_\lambda(u):=\lambda u,$  $u\in E_a,$ $a\in V.$ This map extends to a holomorphic automorphism of $\overline E$. 

Let $V_0$ be an open subset of $V$ which is naturally  identified with an open subset of the section 0 in $E$. A diffeomorphism $\tau$ from a neighbourhood of $V_0$ in $X$ to a neighbourhood of $V_0$ in $E$ is called 
{\it admissible} if it satisfies  essentially the following three conditions:
the restriction of $\tau$ to $V_0$ is the  identity,  the differential  of $\tau$ at each point  $a\in V_0$  is
$\C$-linear and the composition  of 
$$E_a\hookrightarrow \Tan_a(E)\to \Tan_a(X)\to E_a$$    
is the identity, where the morphism $\Tan_a(E)\to \Tan_a(X)$ is given by the differential  of $\tau^{-1}$ at   $a$ and the other maps are the canonical ones, see  \cite{DinhSibony5} for  details. 

Note  that an admissible map is not necessarily holomorphic. 
When $V_0$ is small enough, 
 there are local holomorphic coordinates on a small neighbourhood $U$ of $V_0$ in $X$ so that over $V_0$ we identify naturally $E$ with $V_0\times \C^{k-l}$ and $U$ with an open neighbourhood of $V_0\times\{0\}$  in $V_0\times \C^{k-l}$ (we reduce $U$ if necessary). 
 In this picture, the identity is a holomorphic admissible map.
 There always  exist admissible maps for $V_0:=V.$
 However, such a global  admissible map is rarely  holomorphic.  

Consider an admissible map $\tau$ as above.
Let $T$ be a positive closed $(p,p)$-current on $X$ without mass on $V$ for simplicity.
Define $T_\lambda:=(A_\lambda)_*\tau_*(T)$. 
The family $(T_\lambda)$ is relatively compact on $\pi^{-1}(V_0)$ when $\lambda\to\infty$: we can extract convergent subsequences for $\lambda\to\infty$. The limit currents $R$  are positive closed $(p,p)$-currents without mass on $V$ which are $V$-conic, i.e. 
$(A_\lambda)_* R =R$ for any $\lambda\in\C^*,$ in other words, $R$ is 
invariant by $A_\lambda.$  

Such a current $R$ depends on the choice of $\lambda\to\infty$ but it is independent of the choice of $\tau$. This property gives us a large flexibility to work with admissible maps.
In particular, using global admissible maps, we obtains positive closed $(p,p)$-currents $R$ on $\overline E$. 
It is also known that the cohomology class of $R$ depends on $T$ but does not depend on the choice of $R$. 
This class is denoted by $\kappa^V(T)$ and is called {\it the total tangent class} of $T$ with respect to $V$. The currents $R$ are called {\it tangent currents} of $T$ along $V$.  The mass of $R$ and the norm of $\kappa^V(T)$ is bounded by a constant times the mass of $T$.

Let $-h$ denote the tautological $(1,1)$-class on $\overline E$. Recall that $H^*(\overline E,\C)$ is a free $H^*(V,\C)$-module generated by $1,h,\ldots, h^{k-l}$ (the fibers of $\overline E$ are of dimension $k-l$). So we can write
$$\kappa^V(T)=\sum\limits_{j=\max(0,l-p)}^{\min(l,k-p)} \pi^*(\kappa^V_j(T))\smallsmile h^{j-l+p}$$
where  $\kappa^V_j(T)$ is  a  class in $H^{l-j,l-j}(V,\C)$ with the convention  that $\kappa^V_j(T)=0$  outside of the range
$ \max(0,l-p)\leq j\leq  \min(l,k-p) .$

Let $s$ be the maximal integer such that $\kappa_s^V(T)\not=0$. We call it the {\it tangential h-dimension} of $T$ along $V$. 
The class $\kappa^V_s(T)$ is pseudo-effective, i.e. contains a positive closed current on $V.$
 The  tangential h-dimension of $T$ is  also  equal to the maximal integer $s\geq 0$ such that $R\wedge \pi^*(\omega_V^s)\not=0,$
 where $\omega_V$ is any K\"ahler form on $V.$ 
Moreover, if $T_n$ and $T$ are positive closed $(p,p)$-currents on $X$ such that
 $T_n\to T$,  then $\kappa_j^V(T_n)\to 0$ for $j>s$ and any limit class of $\kappa_s^V(T_n)$ is  pseudo-effective
 and is smaller than or equal to $\kappa_s^V(T)$. 
 
 The following result will allow us to bound the number of isolated periodic points of a meromorphic map. We identify here the cohomology group $H^{2k}(X,\C)$ with $\C$ using the integrals of top degree differential forms on $X$. 
 
 \begin{proposition} \label{prop_intersection}
 Let $\Gamma_n$ be complex subvarieties of pure dimension $k-l$ in $X$. Assume that there is a sequence of positive numbers $d_n$ such that $d_n\to \infty$ and $d_n^{-1}[\Gamma_n]$ converges to a positive closed $(l,l)$-current $T$ on $X$. Assume also that the h-tangent dimension of $T$ with respect to $V$ is $0$ and that $\{T\}\smallsmile \{V\}=1$. Then the number $\delta_n$ of isolated points in the intersection $\Gamma_n\cap V$, counted with multiplicity, satisfies $\delta_n\leq d_n+o(d_n)$ as $n\to\infty$. 
 \end{proposition}
 
 We need the following lemma.
 
 \begin{lemma} \label{lemma_intersection}
 Let $\Gamma$ be a subvariety of pure dimension $k-l$ in $X$. Let $a_1,\ldots, a_N$ be the isolated points in $\Gamma\cap V$ and $m_i$ the  multiplicity of the intersection of $\Gamma\cap V$ at $a_i$. Then any tangent current of $[\Gamma]$ along $V$ is larger than or equal to $\sum m_i[\pi^{-1}(a_i)]$. 
 \end{lemma}
\proof
Consider a small open set $V_0$ in $V$ which contains only one point $a_i$. As above, we identify $E$ (resp. $\overline E$) over $V_0$ with $V_0\times \C^{k-l}$ (resp.  $V_0\times \P^{k-l}$), and a small neighbourhood of $V_0$ in $X$ with an open neighbourhood of $V_0\times \{0\}$ in $V_0\times \C^{k- l},$
and  $\pi$  with the canonical projection of $V_0\times \P^{k-l}$ onto its first factor. The identity is then an admissible map. It is clear in this picture that any tangent current of $[\Gamma]$ along $V$ constructed as above is larger than or equal to $m_i[\pi^*(a_i)]$. The lemma follows. 
\endproof
 
 \noindent
{\bf Proof of Proposition \ref{prop_intersection}.}
Define $T_n:=d_n^{-1}[\Gamma_n]$. Extracting a subsequence allows us to assume that $\kappa^V(T_n)$ converges to a class $\kappa$. Since the h-tangent dimension of $T$ is zero, the above discussion implies that
$\kappa=\lambda c$, where $c$ is the class of a fiber of $\overline E$ and $\lambda$ is a positive number.
We also have $\kappa^V(T)=\pi^*(\kappa^V_0(T))$.   
In the above construction of $\kappa^V(T)$ with a global admissible map, we see that the de Rham cohomology class of $T_\lambda$ in a neighbourhood of $V_0\times\{0\}$ does not depend on $\lambda$ when $\lambda\to\infty$. It follows that  $\{T\}\smallsmile \{V\} =\kappa^V(T)\smallsmile \{V\}$. This together with the hypothesis 
$\{T\}\smallsmile \{V\}=1$ implies that
$\kappa^V(T)=c$. The above discussion on the upper semi-continuity of $\kappa_s^V(T_n)$ implies that $\lambda\leq 1$.

By Lemma \ref{lemma_intersection}, we can write $\kappa^V(T_n)=\delta_n d_n^{-1} c+c_n$ where  $c_n$ is a pseudo-effective class. Since $\kappa^V(T_n)$ converges to $\kappa=\lambda c$, we deduce that the cluster values of $c_n$ are also equal to positive constants times $c$ and then $\limsup \delta_n d_n^{-1}\leq \lambda$. The proposition follows. 
 \hfill $\square$


\section{Construction of good inverse branches}  \label{section_inverse_branches}

Consider a map $f:X\to X$ as above with dominant topological degree. 
The purpose of this section is to construct for generic small balls an almost maximal number of inverse branches with respect to $f^n$ that we control the size.

Recall that $I,I',d_p,d_t,\Gamma$ denote the indeterminacy sets, the dynamical degree of order $p$, the topological degree and the closure of the graph of $f$ in $X\times X$. By definition, the dynamical degree of order $p$ and the topological degree of $f^n$ are equal to $d_p^n$ and $d_t^n$ respectively. 
Denote by $I(f^n),I'(f^n),\Gamma_n$ the indeterminacy sets and the closure of the graph of $f^n$.  
The natural projections from $X\times X$ onto its factors are denoted by $\pi_1$ and $\pi_2$. Recall also that 
$I'_0:=I'$, $I'_{n+1}:=I'_0\cup f(I'_n)$ for $n\geq 0$ and $I'_\infty:=\cup_{n\geq 0} I'_n$. 
One should distinguish 
$I'_\infty$ with the set $\cup_{n\geq 0} f^n(I')$ and the union of $I'(f^n)$ which are a priori smaller. 

Choose an analytic subset $\Sigma_0$ of $X$ containing $I,I',f(I),f^{-1}(I')$ such that $\pi_2$ restricted to $\Gamma\setminus\pi_2^{-1}(\Sigma_0)$ defines an unramified covering over $X\setminus \Sigma_0$. Let $B$ be a connected subset of $X$, e.g.
a holomorphic ball, a holomorphic disc or a family of discs through a point in $X$.  
We call\footnote{We can weaken the conditions in the definition but the ones given here are simple and sufficient for our purpose.}
  {\it an inverse branch of order $n$} of $B$ any continuous bijective map $g:B\to B_{-n}$ with $B_{-n}\subset X$ such that if we define $B_{-i}:=f(B_{-i-1})$ with $0\leq i\leq n-1$, then $B_{-i}\cap \Sigma_0=\varnothing$ for $0\leq i\leq n$,  $f:B_{-i}\to B_{-i+1}$ is a bijective map
for $1\leq i\leq n$, $B_0=B$ and $f^n\circ g=\id$ on $B$.

Note that  $f^{n-i}\circ g:B\to B_{-i}$ is an inverse branch of order $i$ of $B$ and
$B$ admits at most $d_t^n$ inverse branches of order $n$.
The condition $B_{-i}\cap \Sigma_0=\varnothing$ implies that the inverse branch can be extended to any small enough open set containing $B$ using local inverses of the map $f^n$. 
We say that the above inverse branch is {\it of size smaller than $\lambda$} if the diameter of $B_{-n}$ is smaller than $\lambda$. We also call $B_{-n}$ {\it the image} of the inverse branch $g:B\to B_{-n}$.

\begin{proposition} \label{prop_branches}
There is a positive closed $(1,1)$-current $R$ on $X$ satisfying the following property. Let $\epsilon,\nu$ be  strictly positive  numbers with $\nu\leq 1$   and let $a$ be a point in $X$ such that the Lelong number $\nu(R,a)$ of $R$ at $a$ is smaller than $\nu$. Then there is a constant $r>0$ such that the ball $B(a,r)$ of center $a$ and of radius $r$ admits at least $(1-\nu)d_t^n$ inverse branches of order $n$ and of size smaller than $(d_{k-1}/d_t+\epsilon)^{n/2}$ for every $n\geq 0$.
\end{proposition}

A theorem by Siu says that $\{\nu(R,a)\geq c\}$ is a proper analytic subset of $X$ for every $c>0$ \cite{Siu}. So the above proposition applies for generic points $a$ in $X$.
We will see later in the construction of $R$ that the set  $\{\nu(R,a)>0\}$ contains the orbits of the critical values and of the indeterminacy points which are obviously an obstruction to obtain inverse branches of balls. However, $\{\nu(R,a)>0\}$ contains a priori other analytic sets which are a less obvious obstruction to the existence of inverse branches. It can be seen as an accumulation locus of the orbits of the indeterminacy points.
We give now the proof of the above proposition. The first step is to define the current $R$. 

Recall that the operators $(f^n)_*$ act continuously on positive closed $(1,1)$-currents and these actions are compatible with the actions of $(f^n)_*$ on $H^{1,1}(X,\R)$. If $T$ is a positive closed $(1,1)$-current on $X$, its mass depends only on the cohomology class $\{T\}$. Therefore, for a fixed norm on cohomology, the mass of $T$ is comparable with the norm of $\{T\}$. We then deduce the existence of
a constant $c_0>0$ independent of $T,f$ and $n$ such that 
$$\|(f^n)_*(T)\|\leq c_0 \|(f^n)_*\|_{H^{1,1}(X,\R)}\|T\|.$$ 
By definition of $d_{k-1}$, we can fix an integer $N\geq 1$ large enough such that 
$c_0\|(f^N)_*\|_{H^{1,1}(X,\R)}<(\theta d_t)^N$, where $d_{k-1}/d_t<\theta<d_{k-1}/d_t+\epsilon$ is any fixed constant strictly smaller than 1.

Define $\Sigma_{i+1}:=f(\Sigma_i)$ for $0\leq i\leq N-1$ and $\Sigma:=\cup_{0\leq i\leq N} \Sigma_i$.  
So any connected and simply connected set outside $\Sigma$ admits the maximal number $d_t^N$ of inverse branches of order $N$ with images outside $\Sigma_0$. Choose a desingularization $\pi:\widehat\Gamma\to \Gamma$ which is a composition of blow-ups of $\Gamma$ along smooth centers in or over $\pi_1^{-1}(\Sigma)\cap\Gamma$ and $\pi_2^{-1}(\Sigma)\cap\Gamma$. Define $\tau_i:=\pi_i\circ\pi$. We can choose $\pi$ so that $\tau_1^{-1}(\Sigma)$ and $\tau_2^{-1}(\Sigma)$ are of pure codimension 1 in $\widehat\Gamma$. By Blanchard's theorem \cite{Blanchard}, $\widehat\Gamma$ is a compact K\"ahler manifold. Fix a K\"ahler form $\widehat\omega$ on $\widehat\Gamma$ which is larger than $\tau_i^*(\omega)$. We also assume that $\widehat \omega$ is large enough so that each ball of radius 1 in $\widehat\Gamma$ with respect to the metric $\widehat \omega$ is contained in an open set biholomorphic to a ball in $\C^k$.

Denote by $\Sigma'$ and $\Sigma''$ respectively the union of components of codimension 1 and the union of components of codimension $\geq 2$ of $\Sigma$. Define 
$$S_0:=c_1\theta^{-N}d_t^N\big([\Sigma']+(\tau_2)_*(\widehat\omega)\big),\qquad 
S:=\sum_{n=0}^N(f_*)^n(S_0)$$
and
$$R:=8\sum_{m\geq 0} (\theta d_t)^{-mN} (f^N)_*^m(S).$$
Here $c_1\geq \delta_1^{-1}$ is a constant satisfying Lemma \ref{lemma_S} below, and  $\delta_1$ is  
a  constant  whose  exact value  will be determined  right after Lemma \ref{lemma_diameter} below.
By definition of $\theta$, the last current is well-defined. Note that one has to distinguish the operators $(f^N)_*^m$ and $(f^{Nm})_*$. The orbit of $\Sigma$ is the obstruction to construct inverse branches of balls. The following lemma shows that it is visible using the current $R$. 

\begin{lemma} \label{lemma_S}
If $c_1$ is large enough, then for every $a\in\Sigma$ the Lelong number $\nu(S_0,a)$ of $S_0$ at $a$ is larger than $1$. 
\end{lemma}
\proof
The lemma is clear for $a\in\Sigma'$. Consider now a point $a\in\Sigma''\setminus\Sigma'$. 
Since the function $\nu(S_0,a)$ is upper semi-continuous in $a$ with respect to the Zariski topology, it is enough to check that $\nu(S_0,a)$ is positive at generic points $a\in \Sigma''\setminus\Sigma'$. 
We then choose $c_1$ large enough in order to get Lelong numbers larger than 1.
So we can assume that $a$ is a regular point of $\Sigma''\setminus\Sigma'$ and  there is a point $\widehat a\in \tau_2^{-1}(a)$ such that $\tau_2^{-1}(\Sigma'')$ is a hypersurface smooth at $\widehat a$.  

Choose local coordinates $\widehat z=(\widehat z_1,\ldots,\widehat z_k)$ on a neighbourhood of $\widehat a$ such that $\widehat z=0$ at $\widehat a$ and $\tau_2^{-1}(\Sigma'')$ is given by $\widehat z_1=0$. Since $\Sigma''$ has codimension $\geq 2$, we can choose $\widehat z$ so that the hyperplanes $\widehat H_\xi:=\{\widehat z_k=\xi\}$ parallel to $\{\widehat z_k=0\}$ are sent to hypersurfaces, denoted by $H_\xi,$  which contain $\Sigma''$ in a neighbourhood of $a$. 

The average of $[\widehat H_\xi]$ with respect to the Lebesgue measure on $\xi$ is a smooth form $\widehat\Theta$. So it is bounded by a constant times $\widehat\omega$. On the other hand, since $[H_\xi]$ has positive Lelong number at $a$, $(\tau_2)_*(\widehat\Theta)$ has positive Lelong number at $a$. We conclude that $(\tau_2)_*(\widehat\omega)$ has positive Lelong number at $a$. This completes the proof of the lemma.
\endproof

We show that $R$ satisfies Proposition \ref{prop_branches}. Fix a point $a$ in $X$ such that $\nu(R,a)\leq \nu$. Fix also local holomorphic coordinates near $a$. We will first construct inverse branches for flat holomorphic discs through $a$ and then extend these inverse branches to a small ball centered at $a$. 
We will identify a neighbourhood of $a$ to the unit ball in $\C^k$ for simplicity.
The following version of Sibony-Wong's theorem is the tool for this extension.

Let $\B_r$ denote the ball of center 0 and of radius $r$ in $\C^k$. The family $\Fc$ of complex lines through 0 is parametrized by the projective space $\P^{k-1}$ which is endowed with the natural Fubini-Study metric. 
This metric induces a natural probability measure on $\Fc$ that we denote by $\Lc$.  
If $\Delta$ is an element of $\Fc$, denote by $\Delta_r$ its intersection with $\B_r$. 

\begin{proposition}\label{prop_SW}
Let $0<\delta_0\leq 1$ be  a constant. Let $\Fc'\subset \Fc$ be such that $\Lc(\Fc')\geq \delta_0$, and $A$ the intersection of $\Fc'$ with $\B_r$. 
Let $h:A\to\C^l$ be a map which is holomorphic on each $\Delta_r$ for $\Delta\in\Fc'$ and which can be extended holomorphically to a neighbourhood of $0$. Then $h$ can be extended to a holomorphic map from 
 $\B_{\lambda r}$ to $\C^l$, where $0<\lambda \leq 1$ 
is a constant depending on $\delta_0$ but independent of $l$, $\Fc'$ and $r$. Moreover, if the extension is still denoted by $h$ then 
$$\sup_{\B_{\lambda r}}\|h-h(0)\|\leq \sup_A \|h-h(0)\|.$$
In particular, if $\|h-h(0)\|<\rho$ and if $h(A)$ does not intersect a complex hypersurface $Z$ of the ball of center $h(0)$ and radius $\rho$, then $h(\B_{\lambda r})$ satisfies the same property.
 \end{proposition}
 \proof
 When $l=1$ the result is due to Sibony-Wong \cite{SibonyWong}. We easily deduce from their result the 
holomorphic extension of $h$ for any dimension $l$. In order to get the inequality in the proposition, assume $h(0)=0$ for simplicity. Let $z$ be a point in $\B(0,\lambda r)$ we have to show that $\|h(z)\|\leq \sup_A \|h\|$. Composing $h$ with a rotation on $\C^l$ allows to assume that $h(z)=(\|h(z)\|,0,\ldots,0)$. We obtain the desired inequality by using Sibony-Wong's theorem for the first coordinate function of $h$. 
 
 For the last assertion, we can write $Z=\{g=0\}$ where $g$ is a holomorphic function on the ball of center $h(0)$ and of radius $\rho$. Sibony-Wong's theorem applied to $1/g\circ h$ implies that $1/g\circ h$ is holomorphic on $\B_{\lambda r}$. Hence $h(\B_{\lambda r})$ does not intersect $Z$.  The proposition follows.
 \endproof

In order to control the size of holomorphic discs, we need the following lemma, see \cite[Lemma 1.55]{DinhSibony4} for the proof which is valid for any compact complex manifold $Y$. 

\begin{lemma} \label{lemma_diameter}
Let $Y$ be a compact complex manifold endowed with a fixed Hermitian metric. Let $\delta_1>0$ be a constant small enough depending on $Y$. 
Let $g:\Delta_r\rightarrow Y$ be a holomorphic map from the disc of center $0$ and of radius $r$ in $\C$ to $Y$. Assume that the area of $g(\Delta _r)$, counted with multiplicity, 
is smaller than $\delta_1$.
Then for any $\epsilon >0$, there is a constant $0<\lambda <1$ independent of $g$ and $r$ such that the diameter of $g(\Delta _{\lambda r})$ is at most equal to $\epsilon \sqrt{\area(g(\Delta _r))}$.  
\end{lemma}

We will apply it to $Y=\widehat\Gamma$. So from now on  we fix  a constant $\delta_1$ satisfying the last lemma for $\widehat\Gamma$.

Note that the current $R$ constructed above can be seen as the obstruction to the existence of good inverse branches for balls in the spirit of Proposition \ref{prop_branches}. In order to measure the obstruction to the inverse branches of discs through a point $a$ we have to slice this current using complex lines through $a$. We will need the following known technical result, see Lemma 5.52  in \cite{DinhSibony4}. Recall that we identify a neighbourhood of $a$ in $X$ to the unit ball in $\C^k$ for simplicity.

\begin{lemma} \label{lemma_slice}
Let $T$ be a positive closed $(1,1)$-current on $X$. Then for any constant $0<\delta_2 <1$ there is a constant  $r>0$ and a family 
$\Fc'\subset \Fc$, such that $\Lc(\Fc')\geq 1-\delta_2$, and for
every $\Delta \in \Fc'$ the measure $T\wedge [\Delta_r]$ is well-defined and of mass smaller than or equal to $\nu (T,a)+\delta_2 $. Here, $\nu (T,a)$ denotes the Lelong number of $T$ at $a$. 
\end{lemma}

Recall that locally we can write $T=\ddc u$ with $u$ a p.s.h. function. The measure  $T\wedge [\Delta_r]$ is well-defined if $u$ is not identically $-\infty$ on $\Delta_r$. This property holds for $\Lc$-almost every $\Delta$ and we have $T\wedge [\Delta_r]:=\ddc (u[\Delta_r])$.

We are now ready to construct inverse branches for discs through the point $a$ under the hypotheses of Proposition \ref{prop_branches}. 
Fix a value  of the constant $c_1\geq \delta_1^{-1}$ in the definition of $S$ satisfying Lemma \ref{lemma_S}.  
We have the following lemma.

\begin{lemma} \label{lemma_branches}
There is a number $r_0>0$ and a family $\Fc_0\subset \Fc$ with $\Lc(\Fc_0)\geq 1-\nu/2$ satisfying the following property.  For every complex line $\Delta\in\Fc_0$ and for every $n\geq 0$, the disc $\Delta_{r_0}$ admits at least $(1-\nu/2)d_t^n$ inverse branches $g:\Delta_{r_0}\to \Delta_{r_0,-n}$ of order $n$ such that if we define $\Delta_{r_0,-i}:= f^{n-i}(\Delta_{r_0,-n})$ for $0\leq i\leq n$ and $\widehat \Delta_{r_0,-i}:=\tau_1^{-1}(\Delta_{r_0,-i-1})$ for $0\leq i\leq n-1$,
then  $\Delta_{r_0,-i}\cap\Sigma=\varnothing$ for $0\leq i\leq n$ and 
the diameters of $\widehat \Delta_{r_0,-i}$ for $0\leq i\leq n-1$  are smaller than ${1\over 2}\theta^{i/2}$. 
\end{lemma}

Note that since $\Delta_{r_0,-i}\cap\Sigma=\varnothing$ for $0\leq i\leq n$, 
$\tau_1$ defines a biholomorphic map between $\widehat \Delta_{r_0,-i}$ and $\Delta_{r_0,-i-1}$ and 
$\tau_2$ defines a biholomorphic map between $\widehat \Delta_{r_0,-i}$ and $\Delta_{r_0,-i}$. 
Moreover, since $\widehat \omega\geq \tau_1^*(\omega)$  and $\widehat \omega\geq \tau_2^*(\omega)$, the diameter of $\widehat \Delta_{r_0,-i}$ is larger than or equal to the diameters of 
$ \Delta_{r_0,-i-1}$ and of  $ \Delta_{r_0,-i}$. So the diameter of $ \Delta_{r_0,-i}$ is smaller than ${1\over 2}\theta^{i/2}$ for $0\leq i\leq n$. 
\proof
Observe that if $g$ is  an inverse branch of order $n$ satisfying the properties in the lemma then $f\circ g$ is an inverse branch of order $n-1$ satisfying the same properties. So we only have to prove the lemma for $n=mN$ where $m$ is an integer. 

By Lemma \ref{lemma_diameter}, we only need to bound the area of $\widehat \Delta_{r_0,-i}$ by $\theta^{i}/c_1\leq \delta_1$ and then reduce the radius $r_0$ in order to get the diameter control. 
The rest of the lemma is obtained by induction on $m$. 
We will only consider discs $\Delta_r$ through $a$ which are not contained in the orbit of $\Sigma$. This property holds for almost every disc. 

By Lemma \ref{lemma_slice} applied to $R$ and to $\delta_2:=\nu/2$, we can choose a number $r$ and a family $\Fc_0\subset \Fc$ with $\Lc(\Fc_0)\geq 1-\nu/2$  such that  for $\Delta\in\Fc_0$ the measure $R\wedge [\Delta_r]$ is well-defined and of mass smaller than $2\nu$. Let $\nu_m$ denote the mass of $d_t^{-Nm} (f^N)_*^m(S)\wedge [\Delta_r]$. By definition of $R$, we have $\sum_{m\geq 0} \theta^{-Nm} \nu_m\leq \nu/4$. 
 
We show by induction on $m$ that  for every $\Delta\in\Fc_0$ the disc $\Delta_r$ admits at least $\gamma_m:=(1-2\sum_{0\leq i<m} \theta^{-Ni} \nu_i)d_t^{Nm}$ inverse branches $g^{(s)}_{-n}:\Delta_r\to \Delta_{r,-Nm}^{(s)}$  of order $Nm$ such that the area of $\widehat \Delta_{r,-i}^{(s)}$ is smaller than $\theta^i/c_1$ and $\Delta_{r,-i}^{(s)}\cap\Sigma=\varnothing$ for $0\leq i\leq Nm$. 
We used here similar notation as the one introduced in the statement of the lemma. The index $s$ satisfies $1\leq s\leq s_m$ for some integer $s_m$ with $\gamma_m\leq s_m\leq d_t^{Nm}$. Then the above discussion implies the result. Assume this property for $m$. The case $m=0$ is clear since the choice of $\Fc_0$  implies that $\Delta_r$ is out of $\Sigma$. We construct inverse branches of order $N(m+1)$.

The property $\Delta_{r,-Nm}^{(s)}\cap \Sigma=\varnothing$ and the definition of $\Sigma$ allow us to define the maximal number $d_t^{N}$ inverse branches of order $N$ for each $\Delta_{r,-Nm}^{(s)}$. Composing them with the inverse branches of order $Nm$ of $\Delta$ gives $\gamma_md_t^N$ inverse branches of order $N(m+1)$ for $\Delta$. We will count and remove the ones which do not satisfy the area control. 
We call them {\it large-sized inverse branches}. We also have to remove later inverse branches whose images intersect $\Sigma$.
We first count the number of large-sized inverse branches of order $N$ of $\Delta_{r,-Nm}^{(s)}$ for each $s$. 
For simplicity, we will drop the letter $s$ for the moment.

Consider all inverse branches $g:\Delta_{r,-Nm}\to \Delta_{r,-Nm-n}$ of order $1\leq n\leq N$ of $\Delta_{r,-Nm}$ such that the area of $\widehat \Delta_{r,-Nm-i}$ is bounded by  $\theta^{Nm+i}/c_1$ for $i\leq n-1$ but not for $i=n$. They are completely disjoint in the sense that such two different branches are not extensions of the same branch of lower order of $\Delta_{r,-Nm}$. The  extensions of order $N(m+1)$ of these branches are exactly the large-sized ones. 
So the number of large-sized branches of order $N(m+1)$ extending $g$ is $d_t^{N-n}$. 

Observe that the area of $\widehat \Delta_{r,-Nm-n}$ is the mass of $[\Delta_{r,-Nm-n}]\wedge (\tau_2)_*(\widehat\omega)$. This mass is smaller than or equal to the mass of $[\Delta_{r,-Nm}]\wedge (f_*)^n(\tau_2)_*(\widehat\omega)$ because $f^n$ defines a biholomorphic map from a neighbourhood of $\Delta_{r,-Nm-n}$ to a neighbourhood of $\Delta_{r,-Nm}$. So the sum of these areas over all these branches $g$ (there are at most $d_t^N$ such maps) is bounded by $c_1^{-1}\theta^N$ times the mass of $S\wedge [\Delta_{r,-Nm}]$. 
Since the area of $\widehat \Delta_{r,-Nm-n}$ is larger than $\theta^{Nm+n}/c_1$,
 the number of considered maps $g$ is at most equal to 
$\theta^{-Nm}$ times the mass of $S\wedge [\Delta_{r,-Nm}]$. 
The number of large-sized inverse branches of order $N$ of $\Delta_{r,-Nm}$ to remove is at most equal to 
$\theta^{-Nm}d_t^N$ times the mass of $S\wedge [\Delta_{r,-Nm}]$. 

Now, the number $M$ of all large-sized inverse branches of order $N(m+1)$ of $\Delta_r$ to remove is bounded by 
$\theta^{-Nm}d_t^N$ times the mass of $\sum_s S\wedge [\Delta^{(s)}_{r,-Nm}]$. By the  definition of inverse branches, the last mass is bounded by the one of $(f^N)^m_*(S)\wedge [\Delta_{r}]$ which is equal to $d_t^{Nm}\nu_m$. We conclude that $M\leq \theta^{-Nm}d_t^{N(m+1)}\nu_m$.  Therefore, the number of inverse branches of order $N(m+1)$ satisfying the area control is larger than or equal to $\gamma_md_t^N-M\geq \gamma_{m+1}+\nu_{m}d_t^{N(m+1)}$. 

It remains to remove the inverse branches whose images intersect $\Sigma$. 
Denote by $t_{m+1}$ the number of inverse branches $g:\Delta_r\to\Delta_{r,-N(m+1)}$ of order $N(m+1)$
such that $\Delta_{r,-Nm-i}\cap \Sigma\not=\varnothing$ for some $1\leq i\leq N$.
By Lemma \ref{lemma_S}, the intersection of $[\Delta_{r,-Nm-i}]$ with the current $S_0$ is a positive measure of mass at least equal to 1.

By definition of inverse branches, the map $f^n$ is holomorphic and injective on a neighbourhood of 
$\Delta_{r,-n}$ with image in a neighbourhood of $\Delta_r$ for every $n\leq N(m+1)$. We then deduce  that  the mass of 
$(f^N)_*^m(S)\wedge [\Delta_r]$ is at least equal to $t_{m+1}$. It follows that $t_{m+1}\leq \nu_m d_t^{Nm}$. 
We conclude that the number of inverse branches of order $N(m+1)$ satisfying the properties in the lemma is at least equal to $\gamma_{m+1}$. 
This completes the proof of the lemma.
\endproof

\medskip
\noindent
{\bf End of the proof of Proposition \ref{prop_branches}.} We now apply 
Lemma \ref{lemma_branches} and
Proposition \ref{prop_SW} for $\delta_0=\nu/4$. Let $a^{(1)}_{-n},\ldots, a^{(s)}_{-n}$ with $0\leq s\leq d_t^n$ be the distinct points in $f^{-n}(a)$ such that $f^i(a^{(j)}_{-n})\not\in \Sigma$ for all $0\leq i\leq n$ and $1\leq j\leq s$. If $g:\Delta_{r_0}\to \Delta_{r_0,-n}$ is an inverse branch as in Lemma \ref{lemma_branches}, then $\Delta_{r_0,-n}$ contains exactly one of the points $a^{(j)}_{-n}$.  We say that $a^{(j)}_{-n}$ is the center of $\Delta_{r_0,-n}$.

Denote by $\Fc^{(j)}$ the set of $\Delta\in\Fc$ such that $\Delta_{r_0}$ admits an inverse branch of order $n$ as in Lemma \ref{lemma_branches} with center $a^{(j)}_{-n}$. Let $S_n$ denote the set of all $j$ such that $\Lc(\Fc^{(j)})\geq \nu/4$. Let $j$ be an element of $S_n$. We show that $B(a,r)$ admits an inverse branch of order $n$ of size $\leq \theta^{n/2}$ with center $a^{(j)}_{-n}$ for a suitable constant $r>0$ independent of $n$.  

Let $A^{(j)}$ denote the intersection of $\Fc^{(j)}$ with $B(a,r_0)$. The inverse branches of $\Delta_{r_0}$ with $\Delta\in\Fc^{(j)}$ with images centered at $a_{-n}^{(j)}$ agree at the common intersection point $a$ and form a map  $g:A^{(j)}\to A^{(j)}_{-n}$ where $A^{(j)}_{-n}$ is the union of $\Delta_{r_0,-n}$ centered at $a^{(j)}_{-n}$ with $\Delta\in \Fc^{(j)}$. Define as above 
$a^{(j)}_{-i}:=f^{n-i}(a^{(j)}_{-n})$,   $A^{(j)}_{-i}:=f^{n-i}(A^{(j)}_{-n})$ for $0\leq i\leq n$ and
$\widehat a^{(j)}_{-i}:=\tau_1^{-1}(a^{(j)}_{-i-1})$,   $\widehat A^{(j)}_{-i}:=\tau_1^{-1}(A^{(j)}_{-i-1})$ for $0\leq i\leq n-1$.   

By Lemma \ref{lemma_branches}, we have $A^{(j)}_{-i}\cap \Sigma=\varnothing$ for $0\leq i\leq n$. Therefore, the map $\tau_1^{-1}\circ f^{n-i-1}\circ g$ extends holomorphically to a neighbourhood of $A^{(j)}$. Moreover, the image of $A^{(j)}$ is equal to $\widehat A^{(j)}_{-i}$ which is contained in the ball of radius ${1\over 2}\theta^{i/2}\leq 1$ centered at $\widehat a^{(j)}_{-i}$ and does not intersect the hypersurface $\tau_1^{-1}(\Sigma)\cup\tau_2^{-1}(\Sigma)$. Recall that the metric $\widehat\omega$ on $\widehat\Gamma$ is chosen so that any ball of radius 1 is contained in an open set biholomorphic to a ball in $\C^k$. So Proposition \ref{prop_SW} can be applied to this map. According to that
proposition,  
for $r$ small enough (equal to $r_0$ times a constant independent of $n,i,j$), all maps $f^{n-i}\circ g$ and $\tau_1^{-1}\circ f^{n-i-1}\circ g$ extend holomorphically to $B(a,r)$. Moreover, their images, denoted by $B(a,r)^{(j)}_{-i}$ and $\widehat B(a,r)^{(j)}_{-i}$ respectively, 
have diameters smaller than or equal to $\theta^{i/2}$. We also have  
$\widehat B(a,r)^{(j)}_{-i}\cap \tau_1^{-1}(\Sigma)=\varnothing$ and
$\widehat B(a,r)^{(j)}_{-i}\cap \tau_2^{-1}(\Sigma)=\varnothing$ for $0\leq i\leq n-1$. It follows that $B(a,r)^{(j)}_{-i}\cap \Sigma=\varnothing$ for $0\leq i\leq n$. So the extension of $g$ defines an inverse branch of order $n$ and of size $\leq \theta^{n/2}$ on $B(a,r)$. 

It remains now to show that $S_n$ contains at least $(1-\nu)d_t^n$ elements. By Lemma \ref{lemma_branches}, we have $\sum_j \Lc(\Fc^{(j)})\geq d_t^n(1-\nu/2)^2$. Since $\Lc(\Fc^{(j)})\leq \Lc(\Fc)=1$, we deduce that 
the last sum is bounded by $\# S_n +(d_t^n-\#S_n)\nu/4$. It follows that $\# S_n +(d_t^n-\#S_n)\nu/4\geq d_t^n(1-\nu/2)^2$. Hence,  $\#S_n\geq (1-\nu)d_t^n$. This completes the proof of the proposition.
\hfill $\square$


\section{Exceptional set for backward orbits} \label{section_equi_preimages}

In this section, we give the proof of Theorem \ref{th_main_2}. In what follows, we only consider points $a$ outside $I_\infty\cup I'_\infty$. We need the following result that was obtained in \cite{DinhSibony3} in a more general setting.

\begin{lemma} \label{lemma_excep_polar}
There is a pluripolar subset $E$ of $X$ containing $I_\infty'$ such that if $a$ is out of $E$ then 
$\mu^a_n$ converges to $\mu$ as $n$ goes to infinity.
\end{lemma}

For every $a\not\in I_\infty'$, define $\epsilon_a:=\sup \|\mu^a-\mu\|$, where the supremum is taken over all cluster values $\mu^a$ of the sequence $\mu^a_n$.  So we have $\mu^a_n\to\mu$ if and only if $\epsilon_a=0$. 
We deduce from the above lemma and Proposition \ref{prop_branches} the following property.

\begin{lemma} \label{lemma_excep_thin}
Let $a$ be a point out of $I_\infty'$. Then, we have $\epsilon_a\leq 2\nu(R,a)$. In particular, $\mu^a_n\to\mu$ if  $\nu(R,a)=0$.
\end{lemma}
\proof
We have seen that the condition $a\not\in I_\infty'$ is necessary to define $\mu^a_n$. Since we always have $\epsilon_a\leq 2$,
we only need to consider the case where $\nu(R,a)<1$. 
Fix a constant $\nu(R,a)<\nu\leq 1$. Let $B(a,r)$ be a ball as in the conclusion of Proposition \ref{prop_branches}. Choose also a point $b$ in $B(a,r)\setminus E$. Such a  choice is  always possible  since $E$ is pluripolar. Write 
$$f^{-n}(a)=a_{-n}^{(1)},\ldots, a_{-n}^{(d_t^n)}$$
and 
$$f^{-n}(b)=b_{-n}^{(1)},\ldots, b_{-n}^{(d_t^n)}$$
where each points are repeated according to its multiplicity. 

Proposition \ref{prop_branches} implies that we can arrange these points so that the distance between $a_{-n}^{(j)}$ and $b_{-n}^{(j)}$ is smaller than $(d_{k-1}/d_t+\epsilon)^{n/2}$ for at least $(1-\nu)d_t^n$ indices $j$. 
Here $\epsilon>0$ is a fixed constant small enough.
Since $(d_{k-1}/d_t+\epsilon)^{n/2}$ tends to 0, we then deduce that any cluster values of the sequence $\mu_n^a-\mu_n^b$ is a measure of mass at most equal to $2\nu$. The property holds for every $\nu>\nu(R,a)$. Hence,  the result follows from Lemma \ref{lemma_excep_polar} which implies that $\mu^b_n\to\mu$. 
\endproof

The exceptional set in Theorem \ref{th_main_2} is given in the following proposition.

\begin{proposition} \label{prop_excep}
There is a proper analytic subset $\Ec$ of $X$, possibly empty, satisfying the following three conditions: 
(1) no component of $\Ec$ is contained in $I_\infty\cup I_\infty'$; (2) $f^{-1}(\Ec\setminus I')\subset \Ec$; (3) any proper analytic subset of $X$ satisfying (1) and (2) is contained in $\Ec$. Moreover, we have 
$$\Ec=\overline{f^{-1}(\Ec\setminus I')} =\overline{f(\Ec\setminus I)}.$$
\end{proposition}
\proof
Consider the set $Y_0:=\{\nu(R,a)\geq 1\}$. By Siu's theorem, $Y_0$ is a proper analytic subset of $X$ \cite{Siu}. Denote for $n\geq 1$ the analytic set $Y_n$ which is the closure of the set
$$\{z\not\in I_\infty\cup I'_\infty,\quad f^{-i}(z)\in Y \text{ for } 0\leq i\leq n\}.$$
Since the sequence $Y_n$ is decreasing, it is stationary: we have $Y_n=Y_{n+1}$ for $n$ large enough. Denote by $\Ec:=Y_n$ for $n$ large enough. 

It is clear that $\Ec$ satisfies the property (1). 
We have by definition  
$$\overline {f^{-1}(\Ec\setminus (I_\infty\cup I'_\infty))}\subset \Ec.$$
Since $f^{-1}(\Ec\setminus (I_\infty\cup I'_\infty))$ is dense in $f^{-1}(\Ec\setminus I')$, the set $\Ec$ satisfies (2). 
Denote by $\Ec_n$ the closure of $f^{-n}(\Ec\setminus (I_\infty\cup I'_\infty))$. This is a decreasing sequence of analytic sets satisfying the property (1). So it is stationary: we have  $\Ec_n=\Ec_{n+1}$ for $m$ large enough. Since $f(\Ec_{n+1}\setminus I)$ is dense in $\Ec_n$, we deduce from the last identity that 
$\Ec_{n-1}=\Ec_n$ and hence, by induction, 
$\Ec_1=\Ec$. It follows that
$\Ec=\overline{f^{-1}(\Ec\setminus I')}$ which also implies that $\Ec =\overline{f(\Ec\setminus I)}$.

Let $\Ec'$ be a proper analytic subset of $X$ satisfying (1) and (2). We have to show that $\Ec'\subset \Ec$. Property (2) implies that if $a$ is a point in $\Ec'\setminus I_\infty'$ then any cluster values of $\mu^a_n$ is supported by $\Ec'$. Since $\mu$ has no mass on $\Ec'$, we deduce that $\epsilon_a=2$. Lemma \ref{lemma_excep_thin} implies that $a$ is in $Y_0$. So we have $\Ec'\subset Y_0$. Property (2) again, together with the definition of $\Ec$, implies that  $\Ec'\subset \Ec$. This completes  the proof of the proposition.
\endproof

We need the following characterization of the exceptional set  $\Ec$.

\begin{lemma} \label{lemma_excep_char}
Let $Y$ be a proper analytic subset of $X$. Let $a$ be a point in $Y$ which does not belong to $I_\infty\cup I_\infty' $. Let $\lambda_n(a)$ denote the number of backward orbits $a_0,a_{-1},\ldots, a_{-n}$ counted with multiplicity with $a_0=a$, $a_{-i-1}\in f^{-1}(a_{-i})$ for $0\leq i\leq n-1$ and $a_{-i}\in Y$ for $0\leq i\leq n$. If $a$ is not in $\Ec$ then $d_t^{-n}\lambda_n(a)\to 0$ as $n\to\infty$.
\end{lemma}
\proof
Observe that since $a$ is out of $I_\infty\cup I_\infty'$ the same property holds for $a_{-i}$. By definition, the sequence $d_t^{-n}\lambda_n(a)$ is decreasing because each backward orbit of order $n+1$ is one of the $d_t$ extensions of backward orbits of order $n$. Since  the functions $\lambda_n$ are upper-semi-continuous with respect to the induced Zariski topology on $Y\setminus (I_\infty\cup I_\infty')$, the function $\lambda:=\lim d_t^{-n}\lambda_n$ is also upper semi-continuous with respect to this topology. Let $m$ denote the maximal value of $\lambda$ on $Y\setminus( I_\infty\cup I_\infty' \cup \Ec)$. It is enough to show that $m=0$. 

Assume that $m>0$. Denote by $Z^*$ the set of points $a\in Y\setminus (I_\infty\cup I_\infty'\cup \Ec)$ such that $\lambda(a)\geq m$. The closure $Z$ of $Z^*$ is an analytic subset of $Y$. No irreducible component of $Z$ is contained in $\Ec$. Consider a point $a\in Z^*$. The invariance properties of $\Ec$ imply that $f^{-1}(a)\cap\Ec=\varnothing$. Using the definition of $\lambda$ and of $m$, we have
$$m=\lambda(a)=d_t^{-1}\sum_{b\in f^{-1}(a)\cap Y} \lambda(b).$$
Since $\lambda(b)\leq m$ and $\#f^{-1}(a)=d_t$, we deduce that $f^{-1}(a)\subset Z$ and $\lambda(b)=m$ for $b\in f^{-1}(a)$. 
Therefore, $f^{-1}(Z^*)\subset Z$ and
$f^{-1}(Z\setminus I')\subset Z$ since 
$f^{-1}(Z^*)$ is dense in $f^{-1}(Z\setminus I')$. Proposition \ref{prop_excep} implies that $Z\subset \Ec$. This is a contradiction. The lemma follows.
\endproof

\noindent
{\bf End of the proof of Theorem \ref{th_main_2}.} 
Let $a$ be a point out of $I_\infty\cup I_\infty'$. If $a$ is in $\Ec$, it is clear that $\mu^a_n$ does not converge to $\mu$. So assume that $a\not\in\Ec$.   We have to show that $\epsilon_a=0$. 
Fix a constant $\nu>0$. It is enough to prove that $\epsilon_a\leq 4\nu$. Define $Y:=\{\nu(R,\cdot)\geq \nu\}$. By Siu's theorem, this is a proper analytic subset of $X$. By Lemma \ref{lemma_excep_thin}, the case where $a\not\in Y$ is clear. From now on assume that $a\in Y$. 
By Lemma \ref{lemma_excep_char} applied to $Y$, we have  $\lambda_m(a)\leq\nu d_t^m$ for some integer $m$ large enough.

Consider all backward orbits of $a$ of order $l\leq m$ of the form
$$\Oc:=\{a_0,a_{-1},\ldots,a_{-l}\} \quad \text{with } a_0=a, f(a_{-i-1})=a_{-i} \text{ for } 0\leq i\leq l-1$$
such that $a_{-i}\in Y$ for $i\leq l-1$ and $a_{-l}\not\in Y$ unless $l=m$. 
These orbits are counted with multiplicity. Using that $\mu^a_n$ is the probability measure equidistributed on $f^{-n}(a)$, it is not difficult to see that 
$$\mu^a_n=\sum_\Oc d_t^{-l}\mu^{a_{-l}}_{n-l}$$
for every $n\geq m$. 
By considering the masses of the measures in the above identity, we have
$$\sum_\Oc d_t^{-l}=1.$$
We then deduce from the same identity that
$$\epsilon_a\leq \sum_\Oc d_t^{-l}\epsilon_{a_{-l}}.$$

Let $\Sigma,\Sigma'$ denote the sets of $\Oc$ with $a_{-l}\in Y$ (hence $l=m$) and with $a_{-l}\not\in Y$ respectively. By definition of $\lambda_n$, we have 
$$\sum_{\Oc\in\Sigma} d_t^{-l}=d_t^{-m}\lambda_m(a)\leq \nu.$$
On the other hand, we have $\epsilon_{a_{-l}} \leq 2$ for every $\Oc$ and
by Lemma \ref{lemma_excep_thin},  $\epsilon_{a_{-l}}\leq 2\nu$ for $\Oc\in\Sigma'$. It follows that 
$$\epsilon_a \leq \sum_{\Oc\in\Sigma} d_t^{-l}\epsilon_{a_{-l}}+\sum_{\Oc\in\Sigma'} d_t^{-l}\epsilon_{a_{-l}}
\leq 2\sum_{\Oc\in\Sigma} d_t^{-l}+2\nu \sum_{\Oc\in\Sigma'} d_t^{-l}
\leq 2\nu+2\nu=4\nu.$$
This completes the proof of the theorem.
\hfill $\square$

\begin{remark} \rm
Define by induction $\Ec_0:=\Ec$, $\Ec_n:=f(\Ec_{n-1})$ and $\Ec_\infty:=\cup_{n\geq 0} \Ec_n$. If $a$ is a point in $\Ec_\infty\setminus I_\infty'$, then $\mu^a_n$ has positive mass on $\Ec$ for some $n\geq 0$. It follows from the invariance properties of $\Ec$ that $\mu^a_n$ does not converge to $\mu$ since $\mu$ has no mass on $\Ec$. Lemma \ref{lemma_excep_char} still holds for $a\not\in\Ec_\infty\cup I_\infty'$. We can show  that $\mu^a_n\to \mu$ for such points $a$. This property is slightly stronger than Theorem \ref{th_main_2}.
\end{remark}


\section{Periodic points and Lyapounov exponents} \label{section_equi_periodic_points}

In this section, we give the proof of Theorem \ref{th_main_1} and then a lower bound for the Lyapounov exponents of the equilibrium measure.

We call {\it fixed point} of $f$ any point $a$ such that $(a,a)$ belongs to the closure $\Gamma$ of the graph of $f$ in $X\times X$. 
A fixed point $a$ is {\it isolated} if $(a,a)$ is isolated in the intersection of $\Gamma$ with the diagonal $\Delta$ of $X\times X$. {\it The multiplicity} of an isolated periodic point $a$ is the multiplicity of the intersection $\Gamma\cap \Delta$ at $(a,a)$. 
A fixed point $a$ is {\it regular} if it is not an indeterminacy point, that is, $a\not\in I.$ Such a point is
 called {\it repelling} if all the eigenvalues of the differential of $f$ at $a$ have modulus strictly larger than $1$. 
So repelling fixed points are isolated with multiplicity 1. 

{\it Periodic points of period $n$} are fixed points of $f^n$. A periodic point $a$ of order $n$ is {\it regular} if $f^i(a)\not\in I$ for every $i\in\N$. Such a point  is said to be 
{\it repelling} if it is moreover a repelling fixed point  of $f^n.$  We need the following  upper bound for the number of isolated periodic points.

\begin{proposition} \label{prop_bound_Pn}
Let $P_n$ denote the number of isolated periodic points of period $n$ of $f$. Then $\#P_n\leq d_t^n+o(d_t^n)$ as $n$ goes to infinity.
\end{proposition}

We first prove the following property.

\begin{lemma} \label{lemma_graphs}
Let $\Gamma_n$ denote the closure of the graph of $f^n$ in $X\times X$. Then the sequence of positive closed $(k,k)$-currents $d_t^{-n}[\Gamma_n]$ converges to $\pi_1^*(\mu)$ as $n$ goes to infinity. Here, $\pi_1:X\times X\to X$ is the natural projection onto the first factor.
\end{lemma}
\proof
Denote by $z=(z_1,z_2)$ a general point in $X\times X$ with $z_1,z_2\in X$.
Observe that smooth $(k,k)$-forms on $X\times X$ can be written as a finite combination of forms of types
$$\Phi(z):=u(z)\Omega(z_1)\wedge\alpha(z_1)\wedge\Theta(z_2)\wedge\beta(z_2),$$
where $u$ is a smooth function, $\Omega,\Theta$ are smooth positive forms and $\alpha,\beta$ are smooth forms of bidegree $(p,0)$ or $(0,p)$ for some $p\geq 0$. We have to check that 
$$\big\langle d_t^{-n}[\Gamma_n]-\pi_1^*(\mu),\Phi\big\rangle \to 0.$$

\smallskip
\noindent
{\bf Case 1.} Assume that $\Omega$ is a function,  $\Theta$ is of bidegree maximal $(k,k)$ and $p=0$. We may, for simplicity, assume that $\Omega=1$, $\alpha=1$ and $\beta=1$. We have by Fubini's theorem
$$\big\langle d_t^{-n}[\Gamma_n],\Phi\big\rangle = \int_{a\in X} d_t^{-n}\sum_{b\in f^{-1}(a)}u(b,a) \Theta(a) = \int_{a\in X}\langle\mu^a_n,u(\cdot,a)\rangle\Theta(a).$$
Since the measure associated to $\Theta$ has no mass on proper analytic subsets of $X$, by Theorem \ref{th_main_2}, for $\Theta$-almost every $a$, we have $\mu^a_n\to \mu$. So the last sum converges to 
$$\int_{a\in X} \langle \mu,u(\cdot,a)\rangle \Theta(a)=\langle \pi_1^*(\mu),\Phi\rangle.$$
The lemma holds in this case. 

\smallskip
\noindent
{\bf Case 2.} Assume that $\Omega$ is of bidegree $(l,l)$ with $l\geq 1$ and $p=0$. For simplicity, we can assume that $\alpha=1$, $\beta=1$, $|u|\leq 1$, $\Omega(z_1)\leq \omega^l(z_1)$ and $\Theta(z_2)\leq \omega^{k-l}(z_2)$. Since $\langle \pi_1^*(\mu),\Phi\rangle=0$ because of bidegree reason on variable $z_1$, we have to verify that $\langle d_t^{-n}[\Gamma_n],\Phi\rangle\to 0$. We have 
$$\Big|\langle d_t^{-n}[\Gamma_n],\Phi\rangle\Big| \leq \langle d_t^{-n}[\Gamma_n],\omega(z_1)^l\wedge\omega(z_2)^{k-l}\rangle = d_t^{-n} \int_X (f^n)^*(\omega^{k-l})\wedge\omega^l.$$
Clearly, the last integral tends to $0$ since $f$ is with dominant topological degree.

\smallskip
\noindent
{\bf Case 3.} In this last case, assume that $p\geq 1$. We also have $\langle \pi_1^*(\mu),\Phi\rangle=0$ because of bidegree reason on variable $z_1$. We check that $\langle d_t^{-n}[\Gamma_n],\Phi\rangle\to 0$. By Cauchy-Swcharz's inequality, we have 
$$\Big|\langle d_t^{-n}[\Gamma_n],\Phi\rangle\Big|^2 \leq \Big|\big\langle d_t^{-n}[\Gamma_n], \Phi_1\big\rangle \Big| 
\Big|\big\langle d_t^{-n}[\Gamma_n], \Phi_2\big\rangle \Big|$$
with
$$\Phi_1:=|u|^2 \alpha(z_1)\wedge\overline{\alpha(z_1)} \wedge \Omega(z_1)\wedge\Theta(z_2) \quad \text{and} \quad \Phi_2:= \beta(z_2)\wedge\overline{\beta(z_2)} \wedge \Omega(z_1)\wedge\Theta(z_2).$$
Using the previous cases, we see that the first factor in the last product tend to $0$ and the second one is bounded. The lemma follows.
\endproof

\noindent
{\bf End of the proof of Proposition \ref{prop_bound_Pn}.} By Proposition \ref{prop_intersection} and Lemma \ref{lemma_graphs}, it is enough to check that $\{\pi_1^*(\mu)\}\smallsmile \{\Delta\}=1$ and that the h-tangent dimension of $\pi_1^*(\mu)$ with respect to $\Delta$ is 0. Since $\mu$ is a probability measure, its class in $H^{k,k}(X)$ is also the class of a point. So the class of $\pi_1^*(\mu)$ is also the class  of fiber of $\pi_1$. Any fiber of $\pi_1$ intersects $\Delta$ transversally at a point. So $\{\pi_1^*(\mu)\}\smallsmile \{\Delta\}=1$.

To calculate  the h-tangent dimension of $\pi_1^*(\mu) ,$ consider a point $a$ in $X$. We identify a neighbourhood of $a$ with a domain $U$ in $\C^k$ endowed with the standard coordinates $z=(z_1,\ldots, z_k)$. They induce a local coordinate system $(z,z')$ on a neighbourhood of the point $(a,a)$ in $\Delta$. We use a new coordinate system $(z,z'')$ with $z'':=z'-z$. In these coordinates, a neighbourhood of $(a,a)$ is identified to an open subset of $U\times\C^k$, $\Delta$ is given by $\{z''=0\}$ and $\pi_1$ is always the natural projection onto $U$. The normal vector bundle of $\Delta$ is  identified to $U\times \C^k$. The identity map is an admissible local map. So it is not difficult to see that the tangent current of $\pi_1^*(\mu)$ along $\Delta$ is also identified to  $\pi_1^*(\mu)$. The h-dimension of this current is clearly 0. This completes the proof of the proposition. 
\hfill $\square$

\medskip

\noindent
{\bf End of the proof of Theorem \ref{th_main_1}.} 
We can now follow the proof for the case of holomorphic maps presented in \cite{DinhSibony4}. 
For the reader's convenience, we give here the details.

First observe that with the notation as in Theorem \ref{th_main_1}, Proposition \ref{prop_bound_Pn} implies that any cluster value of the sequence $\mu_n:=d_t^{-n} \sum_{a\in Q_n} \delta_a$ is a positive measure of mass at most equal to $1.$ Therefore, it suffices to consider the case where $Q_n$ is the smallest set, i.e. the intersection of $RP_n$ with the support $\supp(\mu)$ of $\mu$. Fix a small constant $\nu>0$. It suffices to prove that any cluster value $\mu'$ of $\mu_n$ satisfies $\mu'\geq (1-4\nu)\mu$. 
Let $B$ be an open subset of $X$. We have to prove that $\mu'(B)\geq (1-4\nu)\mu(B)$. 

Recall that  $\mu$ has no mass on proper analytic subsets of $X$. So it has no mass on $I_\infty\cup I_\infty'$ nor on 
$\{\nu(R,\cdot)>0\}$. In what follows, we only consider balls  whose centers  stay  outside these sets and belong to $\supp(\mu)$, in particular, we have 
$\nu(R,a)=0$ for  such a center $a.$ By Proposition \ref{prop_branches},  any small enough ball centered at $a$ admits at least $(1-\nu)d_t^n$ inverse branches of order $n$ of
diameter $\leq (d_{k-1}/d_t+\epsilon)^{n/2}$. The constant $\epsilon$ is fixed so that 
$d_{k-1}/d_t+\epsilon<1$. 
We only consider such balls.

Since these balls cover a set of full $\mu$ measure, we can find in $B$ a disjoint union of open sets of total $\mu$ measure such that each of these open sets is contained in one of the above considered balls and is biholomorphic to a cube in $\C^k$. Therefore, for simplicity, we can assume that $B$ is such a cube.
We only have to
check for $\mu$-almost every point $a\in X$ that $\#Q_n\cap B\geq
(1-4\nu)d_t^{n}\mu(B)$ when $n$ is large enough.

Choose a finite family of balls $B_i$ of center
$b_i$ with $1\leq i\leq m$
such that $\mu(B_1\cup\ldots\cup B_m)> 1-\nu$ and each
$B_i$ admits $(1-\nu\mu(B))d_t^{n}$ inverse branches of order $n$ with
diameter $\leq (d_{k-1}/d_t+\epsilon)^{n/2}$ for $n$ large enough. Choose
balls $B_i'\Subset B_i$ such that $\mu(B'_1\cup\ldots\cup B'_m)>
1-\nu$. 

Fix a constant $N$ large enough. Since $d_t^{-n} (f^n)^*(\delta_a)$ converge to
$\mu$ for a generic point $a$ in $B$, 
the fiber $f^{-N}(a)$ contains at least $(1-\nu)d_t^N$ points in $\cup B_i'$. Therefore, 
$B$ admits at least $(1-2\nu)d_t^{N}$ inverse branches of order $N$ with small diameters
whose images intersect $\cup B_j'$. So 
the image of such a branch should be contained in one of the $B_j$. 
Choose an open set $B'\subset B$ such that $\mu(B')>(1-\nu)\mu(B)$.
In the same way, we show
that for $n$ large enough,
each $B_j$ admits $(1-2\nu)\mu(B)d_t^{n-N}$ inverse branches 
of order $n-N$ whose images intersect $B'$ and hence are contained  
in $B$. Observe that composing  an inverse branch  of order  $N$ of $B$ whose image 
is contained in $B_j$  with
an inverse  branch of order  $n-N$ of $B_j$  whose image is  contained  in $B,$ we 
obtain  an inverse branch  of order $n$  of $B$  whose image is  contained in $B.$
Consequently, it follows that  $B$ admits at least 
$(1-2\nu)^2\mu(B)d_t^n$ inverse branches $g_i:B\rightarrow
B^{(i)}$ of order $n$ with image $B^{(i)}\Subset B$. 

Every 
holomorphic map $g:U\rightarrow V\Subset U$ on a convex open subset $U$ of $\C^k$ contracts the Kobayashi
metric on $U$ and then admits an attractive fixed point $z$. Moreover, $g^l$
converges uniformly to $z$ and $\cap_{l\geq 0} g^l(U)=\{z\}$. Therefore, each
$g_i$ admits a fixed attractive point $a^{(i)}$. This point is fixed
and repelling for $f^n$. They are different since the $B^{(i)}$ are disjoint. 
Moreover, by definition of inverse branches, the orbit of $a^{(i)}$ does not intersect $\Sigma_0\supset I$. So this is a repelling periodic point of period $n$ for $f$ in our sense.
 
Finally, since $\mu$ is totally invariant, its support satisfies $f^{-1}(\supp(\mu)\setminus I')\subset \supp(\mu)$. Hence, $a^{(i)}$, which is equal to $\cap_{l\geq 0}
g_i^l(\supp(\mu)\cap B)$, is necessarily in
$\supp(\mu)$. 
We deduce that 
$$\#Q_n\cap B\geq (1-2\nu)^2\mu(B)d_t^{n}
\geq (1-4\nu)d_t^{n}\mu(B).$$ 
This completes the proof.
\hfill $\square$

\medskip

\begin{remark} \rm
By Schwarz's lemma, since the diameter of $B^{(i)}$ is smaller than or equal to $(d_{k-1}/d_t+\epsilon)^{n/2}$ all eigenvalues of the differential of $g_i$ at $a^{(i)}$ have modulus smaller than  or equal to this constant. We deduce that the eigenvalues of the differential of $f^n$ at $a^{(i)}$ have modulus larger than or equal to  $(d_{k-1}/d_t+\epsilon)^{-n/2}$. Denote by $Q_n^\epsilon$ the set of repelling periodic points in $Q_n$ satisfying the last property. We then have
$$d_t^{-n} \sum_{a\in Q_n^\epsilon} \delta_a\to \mu.$$
\end{remark}

Using Proposition \ref{prop_branches}, we obtain the following result as in the case of holomorphic maps. 
 
\begin{theorem} \label{th_main_3}
Let $f:X\to X$, $d_t$, $d_p$ and $\mu$ be as in Introduction.   Then  the measure $\mu$ is hyperbolic. Its Lyapounov exponents are bounded from below by ${1\over 2}\log{d_t\over d_{k-1}}$ which is a strictly positive number.
\end{theorem}

The result was  stated for projective manifolds in \cite{Guedj1} but its proof 
 is incomplete since the author uses again his lemma mentioned in the introduction.
It was also stated in \cite{DinhSibony4} for the general case of compact K\"ahler manifolds. 
We give here the details for the reader's convenience, see \cite[Th. 1.120]{DinhSibony4}. Note that the result can be also deduced from a more recent theorem by de Th\'elin \cite{deThelin}.
The cases of endomorphisms of $\P^k$ and of polynomial-like maps were obtained in \cite{BriendDuval1, DinhSibony0}. See also \cite{Lyubich} for the dimension 1 case.

\proof
 Recall that  quasi-p.s.h. functions are $\mu$-integrable. Let $J( f )$ be the Jacobian of $f$ with respect to the K\"ahler metric $\omega$ on $X$. Then, using a resolution of singularity for the graph of $f$ and local holomorphic coordinates, it is not difficult to show that
$|\log J(f)|\leq |\varphi|$ for some quasi-p.s.h. function $\varphi$. So $|\log J(f)|$ is integrable with respect to $\mu$ and
therefore, we can apply Oseledec's
theorem to the natural extension of $f$ (a natural invertible map defined on the space of backward orbits of $f$) \cite{KatokHasselblatt}. We deduce from this result that the smallest Lyapounov exponent of
$\mu$  is equal to
$$\chi := \lim_{n\to \infty}
-{1\over  n} \log\| Df^n(x)^{-1}\|$$
for $\mu$-almost every $x,$ where $Df^n$ denotes the differential of $f^n$.   

Fix  a small constant $\epsilon>0.$ By Proposition  \ref{prop_branches},  there is a ball $B$ of positive $\mu$ measure
which admits at least ${1\over 2} d_t^n $ inverse branches $g_i :\ B\to B^{(i)}_{-n}$ of order $n$ with diameter
$\leq  ({d_{k-1}\over d_t}+\epsilon)^{n/ 2}.$ By Cauchy's formula, if we reduce  slightly the ball $B,$ we can assume that $\| Dg_i\| \leq  A  ({d_{k-1}\over d_t}+\epsilon)^{n/2}$  for
some constant $ A > 0.$ We then deduce 
that $\| (Df^n)^{-1}\| \leq  A  ({d_{k-1}\over d_t}+\epsilon)^{n/2} $ on $B^{(i)}_{-n}.$ 

The union $V_n$ of the $B^{(i)}_{-n}$ is of measure at least equal
to ${1\over 2}\mu (B)$ since $\mu$ is totally invariant. Therefore, by Fatou's lemma,
$$ {1\over 2}
\mu (B) \leq  \limsup_{n\to\infty}
\langle \mu,\textbf{1}_{V_n}\rangle
\leq  \langle \mu,\limsup \textbf{1}_{V_n}\rangle
 = \langle \mu,\textbf{1}_{\limsup V_n}\rangle.
$$
Hence, for $x$ in the  set $K := \limsup V_n$, which has  positive $\mu$ measure, we
have $\| (Df^n)^{-1}\| \leq  A  ({d_{k-1}\over d_t}+\epsilon)^{n/2} $ for infinitely many of $n.$
Hence,   $\chi \geq {1\over 2}
\log ({d_t\over d_{k-1}+\epsilon d_t}).$ We obtain the result by letting $\epsilon$ to $0.$
\endproof

\noindent
T.-C. Dinh, UPMC Univ Paris 06, UMR 7586, Institut de
Math{\'e}matiques de Jussieu, 4 place Jussieu, F-75005 Paris, France.\\ 
DMA, UMR 8553, Ecole Normale Sup\'erieure,
45 rue d'Ulm, 75005 Paris, France.\\
{\tt dinh@math.jussieu.fr}, {\tt http://www.math.jussieu.fr/$\sim$dinh}

\medskip

\noindent
V.-A.  Nguy{\^e}n,
 Math{\'e}matique-B{\^a}timent 425, UMR 8628, Universit{\'e} Paris-Sud,
91405 Orsay, France.\\
  {\tt VietAnh.Nguyen@math.u-psud.fr}, {\tt http://www.math.u-psud.fr/$\sim$vietanh}

\medskip

\noindent
T.-T. Truong,
Department of Mathematics, Syracuse University, Syracuse, NY 13244, USA
{\tt tutruong@syr.edu}

\end{document}